\documentclass[a4paper,12pt]{amsproc}

\usepackage{xcolor} 
\usepackage{mathtools, todonotes}

\usepackage[colorlinks, backref]{hyperref}
\usepackage{amsrefs,amssymb}
\usepackage{enumerate}

\newcommand{\ignore}[1]{}

\DeclareMathOperator{\acl}{acl}
\DeclareMathOperator{\stab}{stab}
\DeclareMathOperator{\Aut}{Aut}
\DeclareMathOperator{\Sym}{Sym}
\DeclareMathOperator{\supp}{supp}

\newcommand{\bA}{{\mathfrak A}}
\newcommand{\bB}{{\mathfrak B}}
\newcommand{\bC}{{\mathfrak C}}

\newcommand{\bX}{{\mathfrak X}}

\DeclareMathOperator{\id}{id}

\theoremstyle{plain}
\newtheorem{theorem}{Theorem}[subsection]
\newtheorem{definition}[theorem]{Definition}

\newtheorem{lemma}[theorem]{Lemma}

\newtheorem{corollary}[theorem]{Corollary}
\newtheorem{conjecture}{Conjecture}

\theoremstyle{definition}
\newtheorem{remark}[theorem]{Remark}

\title[Mixed identities for oligomorphic automorphism groups]{Mixed identities for oligomorphic automorphism groups} 


\author{Manuel Bodirsky}
\address{M.~Bodirsky, TU Dresden, 01062 Dresden, Germany} 
\email{manuel.bodirsky@tu-dresden.de}
\thanks{The first author received funding from the European Union (project POCOCOP, ERC Synergy grant No.~101071674) and the DFG (project FinHom, Grant 467967530). Views and opinions expressed are however those of the author(s) only and do not necessarily reflect those of the European Union or the European Research Council Executive Agency. Neither the European Union nor the granting authority can be held responsible for them.}
\author{Jakob Schneider}
\address{J.~Schneider, TU Dresden, 01062 Dresden, Germany}
\email{jakob.schneider@tu-dresden.de}
\author{Andreas Thom}
\address{A.~Thom, TU Dresden, 01062 Dresden, Germany}
\email{andreas.thom@tu-dresden.de}
\thanks{The third author acknowledges
funding by the Deutsche Forschungsgemeinschaft (SPP 2026).}

\begin{document}
\begin{abstract}
We study mixed identities for oligomorphic automorphism groups of countable relational structures. Our main result gives sufficient conditions for such a group to not admit a mixed identity without particular constants. We study numerous examples and prove in many cases that there cannot be a non-singular mixed identity.
\end{abstract}
\maketitle
\tableofcontents

\section{Introduction}

In this paper, we study so-called \emph{mixed identities} for oligomorphic automorphism groups of \emph{countable relational structures}. This continues the work of \cites{bradfordschneiderthom2023non,bradfordschneiderthom2023length} and we shall use the same notation as there: By a \emph{word with constants} in the group $G$ and in $r$ variables $x_1,\ldots,x_r$, we mean an element
$$
c_0x_{\iota(1)}^{\varepsilon(1)}c_1\cdots c_{l-1}x_{\iota(l)}^{\varepsilon(l)}c_l\in G\ast\mathbf{F}_r,
$$
of the free product, where $\mathbf{F}_r=\langle x_1,\ldots,x_r\rangle$ denotes the free group of rank $r$; $\varepsilon(j)\in\{\pm1\}$ ($j\in\{1,\ldots,l\}$), and $c_0,\ldots,c_l\in G$. We call $w$ \emph{reduced} if $\iota(j)=\iota(j+1)$ and $\varepsilon(j)=-\varepsilon(j+1)$ implies that $c_j \notin\mathbf{Z}(G)$ for $j\in\{1,\ldots,l-1\}$. In this article, all words with constants are reduced. 
The \emph{length} of the (reduced) word $w$ is defined to be  $\lvert w\rvert\coloneqq l$. We call the constants $c_1,\ldots,c_{l-1}$ \emph{intermediate constants}. Moreover, we define the sets of indices
\begin{gather*}
J_0(w)\coloneqq\{j\in\{1,\ldots,l-1\}\mid\iota(j)\neq\iota(j+1)\},\\
J_+(w)\coloneqq\{j\in\{1,\ldots,l-1\}\mid\iota(j)=\iota(j+1)\text{ and }\varepsilon(j)=\varepsilon(j+1)\},\\
J_-(w)\coloneqq\{j\in\{1,\ldots,l-1\}\mid\iota(j)=\iota(j+1)\text{ and }\varepsilon(j)=-\varepsilon(j+1)\}.
\end{gather*}
The elements of $J_-(w)$ are called \emph{critical indices} and the elements $c_j$ with $j\in J_-(w)$ are called \emph{critical constants}. The latter will play a crucial role in the course of the article. We call $w$ \emph{strong} if there are no critical indices. In other words, $w$ is strong if it does not contain a subword of the form $x^{-1} c x$ or of the form $x c x^{-1}$ for a variable $x$ and a constant $c$. (In \cite{tomanov1985generalized} the term \emph{strict} is used; however, the term \emph{strong} is used in the  recent articles \cites{schneiderthom2022word,bradfordschneiderthom2023length}.)  

From $w$ we obtain the corresponding word map $w\colon G^r\to G$ by substitution, which is also denoted by $w$: Write $w(g_1,\ldots,g_r)\in G$ for the image of $w$ under the unique homomorphism $G\ast\mathbf{F}_r\to G$ which fixes $G$ elementwise and maps $x_i\mapsto g_i$ ($i\in\{1,\ldots,r\})$. The word $w$ is called a \emph{mixed identity} for $G$ iff $w\neq 1_{G\ast\mathbf{F}_r}$ and $w(G^r)=\{1_G\}$. Strong mixed identities arise naturally from coset identities on finite index subgroups and deserve special attention.

We consider the \emph{natural augmentation} $\epsilon \colon G\ast\mathbf{F}_r \to \mathbf{F}_r$ that sends the elements of $G$ to $1_{\mathbf{F}_r}$ and fixes $\mathbf{F}_r.$ For a word with constants, we call $\epsilon(w)$ the \emph{content} of $w$. A word is called \emph{singular} if $\epsilon(w)=1_{\mathbf{F}_r}$ and \emph{non-singular} otherwise. Note that every strong word of positive length is non-singular.

\vspace{0.2cm}

While the study of finite groups \cites{bradfordschneiderthom2023non,schneiderthom2024simple}, finite-dimensional unitary groups \cite{klyachkothom2017new}, and algebraic groups \cite{gordeevkunyavskiiplotkin2018wordmaps} shows that non-singular word maps tend to have a larger word image, our understanding of non-singular word maps is limited for a number of naturally occurring groups. In view of a resolution of Kervaire's Conjecture and its generalization \cite{klyachkothom2017new}*{Conjecture~1.1} it is of considerable importance to identify interesting groups so that non-singular word maps are surjective or, as a first step, at least non-trivial. A key example is the automorphism group of $(\mathbb{Q},<).$ It was proved by Adeleke--Holland in \cite{MR1269842} that word maps without constants are surjective for $\Aut(\mathbb{Q},<)$. 

Apart from the case $\Sym(\omega)$ \cite{hullosin2016transitivity} and ${\rm GL}(\mathbb F_q^{\oplus \omega})$ \cite{bradfordschneiderthom2023non}, word maps for words with constants in oligomorphic automorphism groups of relational structures have not been studied much in the literature. However, in \cite{etedadialiabadigaolemaitremelleray2021dense}, Etedadialiabadi, Gao, Le Maitre,  and Melleray studied dense locally finite subgroups of automorphism groups and obtained some results on the non-existence of mixed identities for the automorphism groups of various versions of the Urysohn metric space and related structures.

\vspace{0.2cm}

Our main motivation for this article is the following conjecture:

\begin{conjecture} 
\label{conj:autq} Let $(X,R)$ be an countably infinite relational structure with oligomorphic automorphism group. Every mixed identity for the group $\Aut(X,R)$ is singular.
\end{conjecture}

In Section \ref{sec:examples}, we prove Conjecture~\ref{conj:autq} for the groups $\Sym(\omega)$, ${\rm GL}(\mathbb{F}_q^{\oplus \omega})$, ${\rm Sp}(\mathbb{F}_q^{\oplus \omega})$ and ${\rm O}(\mathbb{F}_q^{\oplus \omega})$ and related examples. For the automorphism group of the Rado graph and, more generally, structures whose age has the free amalgamation property, as well as the random poset and the random permutation, we prove that there do not exist any mixed identities at all. 
More boldly, one could ask whether every  word map for a non-singular word with constants in a simple oligomorphic automorphism group as above is surjective. This holds in a different setting for the compact Lie group ${\rm U}(n)$ and in one variable by the famous result of Gerstenhaber--Rothaus \cites{gerstenhaberrothaus1962solution, klyachkothom2017new} and for $\Aut(\mathbb Q,<)$ for words without constants as mentioned above; but any positive result in this direction is beyond the scope of this article.

Conjecture~\ref{conj:autq} is interesting for the group $(\mathbb Q,<)$, and non-singularity is crucial, since it follows from work of Zarzycki \cite{zarzycki2010limits} that $\Aut(\mathbb{Q},<)$ admits a singular mixed identity. We will prove the conjecture above for $(\mathbb Q,<)$ in Section \ref{sec:autQ} under various additional assumptions on the mixed identity $w$, for example that $\epsilon(w) \not \in [[\mathbf{F}_2,\mathbf{F}_2],[\mathbf{F}_2,\mathbf{F}_2]]$, see Theorem~\ref{thm:mixidaut}, or that $w$ is strong. Another sufficient condition involves additional constraints on the critical constants. Indeed, the conjecture also holds for words with constants if none of the critical constants fixes a non-trivial interval in $\mathbb{Q}$, see Theorem~\ref{thm:main}.

The methods to deal with these cases and variations thereof have extensions to automorphism groups of homogeneous relational structures with finite signature,  or, more generally, to  oligomorphic automorphism groups; see Section \ref{sec:finrel}. Indeed, in Theorem~\ref{thm:main} we identify various natural conditions that ensure that a word with constants is not a mixed identity and explain numerous examples of relational structures, see Section \ref{sec:examples}. 

In parallel, we collect results from the literature and explain how they imply the existence of mixed identities for particular automorphism groups, such as the ones of
the generic cyclic order, the generic semi-linear order and the Wa\.{z}ewski dendrite. Many more examples can be treated with the same methods, but we leave it there. 

\section{General homogeneous relational structures}
\label{sec:finrel}
In this section we first recall basic concepts from model theory, then present the main theorem about non-existence of certain mixed identities.

\subsection{Preliminaries}
\label{sect:prelims}
All definitions in this section are standard in model theory (see, e.g.,~\cite{hodges1997shorter}), except for Definition~\ref{def:alg-conv}, which appears to be new. 
Let $(X,R)$ be a pair consisting of a countable set $X$ and a set $R=(R_i)_{i\in I}$ of relations on $X$, with $R_i\subseteq X^{m_i}$ where $m_i\in\mathbb{N}$, for $i\in I$, is called the \emph{arity} of $R_i$. 
We call such a pair $(X,R)$ a \emph{relational structure}, $X$ its \emph{domain}, and $(m_i)_{i \in I}$ its \emph{signature}. 
We say that a relational structure $(X,R)$ is \emph{finite} if its domain $X$ is finite, and that it has \emph{finite signature} if $I$ is finite. 
Subsequently, we will consider the automorphism group $\Aut(X,R)$, which is the set of bijections $g$ of $X$ such that both $g$ and $g^{-1}$ preserve all the relations of $R$, i.e.,
for all $i\in I$ we have that $(x_1,\ldots,x_{m_i}) \in R_i$ holds if and only if $(x_1.g,\ldots,x_{m_i}.g) \in R_i$.

For $Y\subseteq X$, we also consider the pair $(Y,(\left.R\right\rvert_Y)_{i \in I})$ and call it a \emph{substructure} of $(X,R)$. We define the stabilizer of $Y$ by $$\stab(Y)\coloneqq\{g\in\Aut(X,R)\mid y.g=y\text{ for all }y\in Y\}.$$

\vspace{0.2cm}

We start out with some definitions and lemmas, which will be necessary to state and prove our main results.

\begin{definition}
 A relational structure $(X,R)$ is called \emph{homogeneous} if every isomorphism between finite  substructures can be extended to an automorphism of $(X,R)$.
\end{definition}

\begin{remark} A subgroup of ${\rm Sym}(X)$ (which is equipped with the product topology where $X$ is taken to be discrete) is closed if and only if it is the automorphism group of a homogeneous relational structure. 
\end{remark}

Every countable homogeneous relational structure is (up to isomorphism) uniquely given by its \emph{age}, i.e., by the class $\mathcal C$ of all finite structures that are isomorphic to a substructure of $(X,R)$. Then, $(X,R)$ is called the \emph{Fra\"iss\'e limit} of ${\mathcal C}$.
By Fra\"iss\'e's theorem, a class $\mathcal C$ of finite structures of the same signature $(m_i)_{i \in I}$ is the age of a countable homogeneous structure if and only if ${\mathcal C}$ is closed under isomorphism and  substructures, 
has countably many isomorphism classes, and has the so-called \emph{amalgamation property}:
for all $(X_1,R_1), (X_2,R_2) \in {\mathcal C}$ 
such that the identity on $X_0 \coloneqq X_1 \cap X_2$ is an isomorphism between the finite  substructures of $(X_1,R_1)$ and $(X_2,R_2)$ with domain $X_0$, respectively, 
there exists $(Y,S) \in {\mathcal C}$ and for every $i \in \{1,2\}$ and isomorphisms $f_i$ 
between $(X_i,R_i)$ and a substructure of $(Y,S)$ such that
$f_1(x)=f_2(x)$ for all $x \in X_0$.
The structure $(Y,S)$ is then called an \emph{amalgam}. 
Some stronger forms of the amalgamation property will be relevant when discussing concrete examples in later sections. 
For the \emph{strong amalgamation property} we additionally require that $f_1(X_1) \cap f_2(X_2) = f_1(X_0) = f_2(X_0)$. 
For the \emph{free amalgamation property} we even require that the structure $(X_1 \cup X_2,R_1,R_2)$ is in ${\mathcal C}$.
Clearly, the free amalgamation property implies the strong amalgamation property.

\begin{definition}
The permutation group $\Aut(X,R)$ is called \emph{oligomorphic} if 
 the component\-wise action of $\Aut(X,R)$ 
 on  $n$-tuples  has for every $n \in \mathbb{N}$ only finitely many orbits.  
\end{definition}

Throughout the text \emph{definable} will always mean \emph{first-order definable with parameters} and if the parameters are restricted to some subset of the domain it will be mentioned explicitly.

\begin{remark}\label{rem:homogen} Note that the automorphism group of every  homogenenous relational structure with finite signature is oligomorphic. {Conversely, if $\Aut(X,R)$ is oligomorphic, and we add all relations to $(X,R)$ that are first-order definable without parameters in $(X,R)$, then the resulting expanded structure is homogeneous and has the same automorphism group as the original structure $(X,R)$.} 
\end{remark}

\begin{remark}
 A relational structure has an oligomorphic automorphism group if and only if it is \emph{$\omega$-categorical}, i.e., if its first-order theory has at most one countable model up to isomorphism (a theorem of Engeler, Svenonius, and Ryll-Nardzewski; see, e.g.,\ \cite{hodges1997shorter}). 
 Moreover, a relation $S$ is definable in an $\omega$-categorical structure $(X,R)$ without parameters if and only if $S$ is preserved by $\Aut(X,R)$. 
 We do not make real use of the logic perspective 
 here, because with the facts mentioned above all the arguments can be translated into the language of permutation groups.
\end{remark}

The following notions are standard in model theory. 

\begin{definition}
Let $(X,R)$ be a relational structure with an oligomorphic automorphism group and let $Y$ be a subset of $X$.
We set 
 \begin{gather*}
 {\rm acl}(Y)\coloneqq\{x\in X\mid\text{there exists a finite set }W \subseteq Y\\ \text{ such that } x.\stab(W)\text{ is finite}\} 
 \end{gather*}
 and call it the \emph{algebraic closure}.
\end{definition}

For $Y$ finite, we have 
${\rm acl}(Y)=\{x\in X\mid x.\stab(Y) \mbox{ is finite}\}$.
One verifies that ${\rm acl}$ is a closure operator. Note the following simple observation.

\begin{lemma} \label{lem:acl_fin}
 Let $(X,R)$ be a relational structure with oligomorphic automorphism group. If $Y\subseteq X$ is finite, then ${\rm acl}(Y)$ is also finite.
\end{lemma}

\begin{proof}
Indeed, for $n=\lvert Y\rvert$, $\Aut(X,R)$ acts with finitely many orbits on $X^{n+1}$. Hence the stabilizer $\stab(Y)\leq\Aut(X,R)$ also has at most finitely many finite orbits on the elements of $X\setminus Y$. Hence ${\rm acl}(Y)$ is finite.
\end{proof}

\begin{definition}
We say that $c\in\Aut(X,R)$ is \emph{small} if $c\in\stab(Z)$ for some  infinite subset $Z$ of $X$ which is definable
in $(X,R)$. 
\end{definition}

The notion of small elements will play a crucial role in Theorem~\ref{thm:main} below.

\begin{lemma}\label{lem:no_smll_elts}
 Let $(X,R)$ be a homogeneous relational structure with oligomorphic automorphism group. The following are equivalent:
 \begin{enumerate}[$(i)$]
 \item If $Y \subseteq X$ is a finite set, $x \in X \setminus{\rm acl}(Y)$, and $Z \subseteq x.\stab(Y)$ is a finite subset, then $\stab(x.\stab(Y) \setminus Z)$ is trivial.
 \item If $Y\subseteq X$ is a finite set, $x\in X\setminus{\rm acl}(Y)$, then $\stab(x.\stab(Y))$ is trivial.
 \item No non-trivial element of $\Aut(X,R)$ is small.
 \end{enumerate}
\end{lemma}

\begin{proof} 
We start by proving the implication (iii)$\Rightarrow$(i). Since $(X,R)$ is homogeneous and $\Aut(X,R)$ is oligomorphic, the set $x.\stab(Y) \setminus Z$ is infinite and definable with parameters $Y \cup Z$. Thus, if (iii) holds, we conclude that $\stab(x.\stab(Y) \setminus Z)$ is trivial.
The implication (i)$\Rightarrow$(ii) is trivial. 
 
 Let's now prove (ii)$\Rightarrow$(iii).
{Suppose that $W$ is an infinite set which is definable in $(X,R)$ by a formula with parameters from a finite set $Y$. Since $W$ is infinite and ${\rm acl}(Y)$ is finite by Lemma~\ref{lem:acl_fin}, 
 there exists  $x\in W\setminus {\rm acl}(Y)$.} Then $x.\stab(Y) \subseteq W$ is infinite. Now, $\stab(x.\stab(Y))$ being trivial implies that $\stab(W)$ is trivial. Hence, we conclude that (iii) holds.
\end{proof}

\begin{definition}
We say that $(X,R)$ has \emph{no algebraicity} if for every finite subset $Y$ of $X$, we have ${\rm acl}(Y)=Y$, i.e., all finite subsets of $X$ are \emph{algebraically closed}.
\end{definition}

\begin{remark}\label{rem:sap}
It is well known that a homogeneous structure has no algebraicity if and only if its age has the strong amalgamation property (see, e.g.,~\cite{cameron1990oligomorphic}).
\end{remark}



The following definition describes an algebraic property a relational structure can have which allows us to prove very strong results on the non-existence of mixed identities for the corresponding automorphism group.

\begin{definition}\label{def:alg-conv}
 We say that a relational structure $(X,R)$ is \emph{algebraically convex} if for every finite set $Y \subseteq X$, the following implication holds: if $x \in {\rm acl}({\rm acl}(Y \cup \{x\}) \setminus \{x\}),$ then $x \in {\rm acl}(Y)$.
\end{definition}

The name of this property is inspired by the fact that the \emph{convex hull} as a closure operator has the same property. Indeed, removing an extreme point from a convex set preserves convexity. It is clear that no algebraicity implies algebraical convexity. Further examples will be discussed in Sections~\ref{sec:semilinear} and~\ref{sec:dendrite}. 
Examples which are not algebraically convex can be found in Section~\ref{sec:eq} and Section~\ref{sec:vectspace}. 

\begin{definition}
 We say that $c\in\Aut(X,R)$ is \emph{slender} if there exists a finite subset $Y$ of $X$ and an infinite definable set $Z \subseteq X$, such that 
 $$ Z \subseteq \{ x\in X \mid x.c\in{\rm acl}(Y\cup\{x\})\}.$$
\end{definition}

\begin{remark} Every small automorphism is slender. Indeed, if $c \in \stab(Z)$ for some infinite definable set $Z$, then the definition of slender is satisfied for $Y= \varnothing$ and $Z$. The converse holds if $(X,R)$ has no algebraicity.

Indeed, assume that there are a finite set $Y\subseteq X$ and an infinite definable set $Z\subseteq X$ such that 
\begin{align*}
Z & \subseteq \{x\in X \mid x.c \in {\rm acl}(Y\cup\{x\})\} \\
& = \{x\in X \mid x.c \in Y \cup \{x\} \} = Y.c^{-1} \cup \{x\in X \mid x.c=x \}.
\end{align*}
Then 
$c \in \stab(Z\setminus Y.c^{-1})$ and $Z\setminus Y.c^{-1}$ is an infinite definable set, proving that $c$ is small.

Hence, we can think of being slender as a generalization of smallness. Indeed, we think that the definition of slender assures in some way that a set of `generalized fixed points'
is large in the sense that it contains an infinite definable subset.
The concept of a slender automorphism will prove to be useful beyond the algebraically convex setting, see Section~\ref{sec:vectspace}.
\end{remark}

\begin{definition} Let $(X,R)$ be a relational structure. We say that a tuple $(x_1,\dots,x_n) \in X^n$ is \emph{staggered algebraically independent} if $x_i\notin{\rm acl}(\{x_1,\dots,x_{i-1}\})$ for every $i \in \{1,\dots,n\}$. In particular, $x_1 \notin {\rm acl}(\varnothing).$
\end{definition}

Staggered algebraically independent tuples have a concrete meaning in many circumstances. For example, a sequence of vectors $(v_1,\dots,v_n)$ in $\mathbb{F}_q^{\oplus\omega}$ (see Section \ref{sec:vectspace}) is staggered algebraically independent if and only if it is linearly independent.

\subsection{The main theorem}
\label{sect:main}

At this point we are prepared to state and prove our main result.

\begin{theorem} \label{thm:main}
Suppose that $(X,R)$ is a homogeneous relational structure with oligomorphic automorphism group. Let $w \in \Aut(X,R) \ast \mathbf{F}_r$ be a word with constants of positive length. Assume that one of the following conditions holds true:
\begin{enumerate}[$(i)$]
\item $w$ is strong,
\item $w$ has no non-trivial 
small critical constants and $(X,R)$ is algebraically convex,
\item $w$ has no non-trivial slender critical constants.
\end{enumerate}
Then, for every $n \in \mathbb{N}$, there exists a staggered algebraically independent tuple $(\alpha_1,\beta_1$, \dots, $\alpha_n,\beta_n) \in X^{2n}$ and $h_1,\dots,h_r \in \Aut(X,R)$ such that
$$\alpha_i.w(h_1,\dots,h_r) = \beta_i, \quad \forall i \in \{1,\dots,n\}.$$
In particular, $w$ is not a mixed identity.
\end{theorem}

\begin{proof}
We prove (i), (ii), and (iii) simultaneously. 
Since the property of mapping some staggered algebraically independent tuple in the required way is invariant under conjugation, we may assume that the first constant $c_0$ of $w$ is trivial. Thus, without loss of generality, $w$ can be written as
$$w=x_{\iota(1)}^{\varepsilon(1)}c_1 x_{\iota(2)}^{\varepsilon(2)}\cdots x_{\iota(l-1)}^{\varepsilon(l-1)}c_{l-1} x_{\iota(l)}^{\varepsilon(l)}c_l\in\Aut(X,R)\ast\mathbf{F}_r.
 $$
We want to use induction on $n$ to establish the claim. Subsequently, assume we are in the $n$\textsuperscript{th} step of the induction for some $n\geq1$. Let $\Lambda_k^\varepsilon$ be some finite sets specified later which contain elements from the previous steps (if any) of the induction ($\varepsilon\in\{\pm1\}$, $k\in\{1,\ldots,r\}$). In the case $n=1$, we set $\Lambda_k^\varepsilon\coloneqq\varnothing$ for all $\varepsilon$ and $k$. Let $\alpha_n\in X\setminus({\rm acl}(\Lambda_{\iota(1)}^{\varepsilon(1)})\cup{\rm acl}(\{\alpha_1,\beta_1,\ldots,\alpha_{n-1},\beta_{n-1}\}))$ be arbitrary and $\beta_n\in X \setminus {\rm acl}(\{\alpha_1,\beta_1,\ldots,\alpha_{n-1},\beta_{n-1},\alpha_n\})$ to be defined later. For $j\in\{1,\ldots,l\}$ we will choose elements $\omega_{j}^{-\varepsilon(j)}\in X$ with the aim to define maps $h_1,\ldots,h_r\in \Aut(X,R)$ such that\begin{gather*}
 \alpha_n\eqqcolon\omega_{1}^{\varepsilon(1)}\stackrel{h^{\varepsilon(1)}_{\iota(1)}}{\mapsto}
 \omega_{1}^{-\varepsilon(1)}\stackrel{c_1}{\mapsto} \omega_{2}^{\varepsilon(2)}\stackrel{h^{\varepsilon(2)}_{\iota(2)}}{\mapsto}
 \cdots
\stackrel{c_{l-1}}{\mapsto}
\omega_{l}^{\varepsilon(l)}\stackrel{h^{\varepsilon(l)}_{\iota(l)}}{\mapsto}\omega_{l}^{-\varepsilon(l)} \stackrel{c_l}{\mapsto}\beta_n.
\end{gather*}
We see that the elements $\omega_j^{\varepsilon(j)}$ are determined by $\omega_1^{\varepsilon(1)}\coloneqq\alpha_n$ and the equations $\omega_{j+1}^{\varepsilon(j+1)} \coloneqq \omega_j^{-\varepsilon(j)}.c_j$ for $j\in\{1,\ldots,l-1\}$ (and the choice $\omega_j^{\varepsilon(j)}.h_{\iota(j)}^{\varepsilon(j)}=\omega_j^{-\varepsilon(j)}$ for all $j$). 
For $k\in\{1,\ldots,r\}$ and $\varepsilon \in \{\pm 1\}$, define 
\begin{align*} \Omega^{\varepsilon}_k & \coloneqq\{\omega^\varepsilon_j\mid j\in\{1,\ldots,l\}, \iota(j)=k\}\cup\Lambda_k^\varepsilon. 
\end{align*}
For 
$j \in \{0,1,\ldots,l\}$,  define 
$$\Omega_{j,k}^{\varepsilon}\coloneqq\{\omega_{j'}^\varepsilon\mid j'\leq j,\iota(j')=k\}\cup\Lambda_k^\varepsilon.$$ 
Note that here $\Omega_{0,k}^\varepsilon\coloneqq\Lambda_k^\varepsilon$ for all $k$. 
For better readability, in exponents we write in the following $+$ instead of $+1$ and $-$ instead of $-1$.

In order for the $h_k$ to exist, we want to achieve in the course of the inductive process that
\begin{itemize}
    \item the substructures induced by $\Omega^{+}_k$ and $\Omega^{-}_k$ are isomorphic (for all $k\in\{1,\ldots,r\}$) via the natural map sending $\omega^{+}_j$ to $\omega^{-}_j$ for $j\in\{1,\ldots,l\}$ with $\iota(j)=k$, and mapping the set $\Lambda_k^+$ bijectively onto $\Lambda_k^-$;
    \item $\omega_j^\varepsilon\neq\omega_{j'}^\varepsilon$ and $\omega_j^\varepsilon\notin\Lambda_{\iota(j)}^\varepsilon$, for $\varepsilon\in\{\pm1\}$ and $j,j'\in\{1,\ldots,l\}$ such that $j\neq j'$ and $\iota(j)=\iota(j')$. 

\end{itemize}
The proof proceeds by induction on $j\in\{1,\ldots,l\}$. In the $j$\textsuperscript{th} step we extend the isomorphism between the  substructures on $\Omega_{j-1,k}^+$ and $\Omega_{j-1,k}^-$ from the previous step to an isomorphism between the  substructures on $\Omega_{j,k}^+$ and 
$\Omega_{j,k}^-$ (for all $k\in\{1,\ldots,r\}$). For $k\neq\iota(j)$, there is nothing to do and we just define $\Omega_{j,k}^\varepsilon\coloneqq\Omega_{j-1,k}^\varepsilon$ for $\varepsilon\in\{\pm1\}$. If $k=\iota(j)$, we will add the arrow $\omega_j^+\mapsto\omega_j^-$ to our partial isomorphism from the previous step and set $\Omega_{j,k}^\varepsilon=\Omega_{j,\iota(j)}^\varepsilon\coloneqq\Omega_{j-1,\iota(j)}^\varepsilon\cup\{\omega_j^\varepsilon\}$ for $\varepsilon\in\{\pm1\}$. During the inductive process, we want to ensure that $\omega_j^\varepsilon\notin{\rm acl}(\Omega_{j-1,\iota(j)}^\varepsilon)$ for any $\varepsilon$.
For $j=1$, we have by definition
$$
\omega^{\varepsilon(1)}_1\coloneqq\alpha_n\notin{\rm acl}(\Omega_{0,\iota(1)}^{\varepsilon(1)})\cup{\rm acl}(\{\alpha_1,\beta_1,\ldots,\alpha_{n-1},\beta_{n-1}\}).
$$

For the inductive step, we assume that
$$\omega_1^+,\omega_1^-,\ldots,\omega_{j-1}^+,\omega_{j-1}^-,\omega_j^{\varepsilon(j)}
$$ 
are already defined and we want to determine a suitable choice for $\omega_j^{-\varepsilon(j)}$. 
By the inductive assumption, $\omega_j^{\varepsilon(j)}\notin{\rm acl}(\Omega_{j-1,\iota(j)}^{\varepsilon(j)})$ and there is an isomorphism between the  substructures on $\Omega_{j-1,\iota(j)}^{\varepsilon(j)}$ and 
$\Omega_{j-1,\iota(j)}^{-\varepsilon(j)}$. Let $g$ be an extension of this partial isomorphism to an automorphism of $(X,R).$ A priori, the set of valid choices for $\omega_j^{-\varepsilon(j)}$ lies in the definable set 
$$
\Sigma_{j}\coloneqq\omega_{j}^{\varepsilon(j)}.g\stab(\Omega_{j-1,\iota(j)}^{-\varepsilon(j)}) .
$$ 
Since $\omega_{j}^{\varepsilon(j)}.g\notin{\rm acl}(\Omega_{j-1,\iota(j)}^{\varepsilon(j)}).g={\rm acl}(\Omega_{j-1,\iota(j)}^{-\varepsilon(j)})$, the set $\Sigma_j$ is infinite.

Apart from this, we need that if $j\in\{1,\ldots,l-1\}$, then $\omega_{j+1}^{\varepsilon(j+1)}=\omega_j^{-\varepsilon(j)}.c_j\notin{\rm acl}(\Omega_{j,\iota(j+1)}^{\varepsilon(j+1)})$. In order to achieve this, we distinguish the following three cases.
\vspace{0.2cm}

\emph{Case~1: $j\in J_0(w)$.} In this case, we have $\iota(j)\neq\iota(j+1)$ and hence $\Omega_{j,\iota(j+1)}^{\varepsilon(j+1)}=\Omega_{j-1,\iota(j+1)}^{\varepsilon(j+1)}$ which is already fixed. Then we must require that

$$\omega_j^{-\varepsilon(j)}\in \Sigma_{j}\setminus ({\rm acl}(\Omega_{j-1,\iota(j+1)}^{\varepsilon(j+1)}).c_j^{-1}
),$$ 
which is clearly possible since $\Sigma_{j}$ is infinite and ${\rm acl}(\Omega_{j-1,\iota(j+1)}^{\varepsilon(j+1)})$ is finite. 
\vspace{0.2cm}

\emph{Case~2: $j\in J_+(w)$.} 
Then $\iota(j)=\iota(j+1)$ and $\varepsilon(j)=\varepsilon(j+1)$, so that $\Omega_{j,\iota(j+1)}^{\varepsilon(j+1)}=\Omega_{j,\iota(j)}^{\varepsilon(j)}$ which is also already fixed. Then we must have that
$$
\omega_j^{-\varepsilon(j)}\in \Sigma_j\setminus ({\rm acl}(\Omega_{j,\iota(j)}^{\varepsilon(j)}).c_j^{-1}),
$$ 
which is again possible since $\Sigma_j$ is infinite and ${\rm acl}(\Omega_{j,\iota(j)}^{\varepsilon(j)})$ is finite.
\vspace{0.2cm}

Thus, the proof of (i) is complete, since 
Case~1 and Case~2 are the only ones that occur when $w$ is strong as $J_-(w)=\varnothing$.
Indeed, note that in step $l$ (the final one) we have infinitely many choices for $\omega_{l}^{-\varepsilon(l)}$. In particular, we can choose $\omega_l^{-\varepsilon(l)}$ such that $\omega_l^{-\varepsilon(l)}.c_l=\beta_n\notin{\rm acl}(\{\alpha_1,\beta_1,\ldots,\alpha_n\})$ which is a finite set by Lemma~\ref{lem:acl_fin}.
As discussed before, by homogeneity, we obtain suitable $h_1,\ldots,h_r \in \Aut(X,R)$
with $\alpha_n.w(h_1,\ldots,h_{r})=\beta_n.$ 

To finish the proof of (ii), we also need to consider the case $J_-(w)\neq\varnothing$:
\vspace{0.2cm}

\emph{Case~3: We have $j\in J_-(w)$.} Then $\iota(j)=\iota(j+1)$ and $\varepsilon(j)=-\varepsilon(j+1)$. In this case we obtain 
$$
\Omega_{j,\iota(j+1)}^{\varepsilon(j+1)}=\Omega_{j,\iota(j)}^{-\varepsilon(j)}=\Omega_{j-1,\iota(j)}^{-\varepsilon(j)}\cup\{\omega_j^{-\varepsilon(j)}\}.
$$ 
We need to find $\omega_j^{-\varepsilon(j)}$ satisfying
$$
\omega_j^{-\varepsilon(j)}\in\Sigma_j \quad \mbox{and} \quad \omega_j^{-\varepsilon(j)}.c_j \not \in{\rm acl}(\Omega_{j-1,\iota(j)}^{-\varepsilon(j)} \cup \{\omega_j^{-\varepsilon(j)}\}).
$$

Let's fix $\omega\in\Sigma_j$ and recall that by definition $\omega\notin{\rm acl}(\Omega_{j-1,\iota(j)}^{-\varepsilon(j)})$. 
By the assumption that $(X,R)$ is algebraically convex we have
$$
\omega \not \in {\rm acl}\left({\rm acl}(\Omega_{j-1,\iota(j)}^{-\varepsilon(j)} \cup \{\omega\}) \setminus \{\omega\} \right).$$
Consider now 
$$\Sigma_j'\coloneqq \omega.\stab \left({\rm acl}(\Omega_{j-1,\iota(j)}^{-\varepsilon(j)} \cup \{\omega\}) \setminus \{\omega\}\right) \subseteq \Sigma_j,
$$
which is definable and infinite by assumption on $\omega$.
For all $\omega'\in\Sigma'_j,$ we have
$$
{\rm acl}(\Omega_{j-1,\iota(j)}^{-\varepsilon(j)} \cup \{\omega'\}) \setminus \{\omega'\} = {\rm acl}(\Omega_{j-1,\iota(j)}^{-\varepsilon(j)} \cup \{\omega\}) \setminus \{\omega\},
$$
since $\omega'=\omega.s$ for some $s \in \stab({\rm acl}(\Omega_{j-1,\iota(j)}^{-\varepsilon(j)} \cup \{\omega\}) \setminus \{\omega\})$ (so that $s$ fixes 
${\rm acl}(\Omega_{j-1,\iota(j)}^{-\varepsilon(j)} \cup \{\omega\}) \setminus \{\omega\})$ 
pointwise, while it also maps $\omega$ to $\omega'$). 
Since $\Sigma_j'$ is infinite, the definable set
$$\Sigma'_j \setminus \left( {\rm acl}(\Omega_{j-1,\iota(j)}^{-\varepsilon(j)} \cup \{\omega\}) \setminus \{\omega\}\right).c_j^{-1}$$
is also infinite.
Since $c_j$ is not small,
we find $$\omega_j^{-\varepsilon(j)} \in \Sigma_j'\setminus \left( {\rm acl}(\Omega_{j-1,\iota(j)}^{-\varepsilon(j)} \cup \{\omega\}) \setminus \{\omega\}\right).c_j^{-1}$$ which is not fixed by $c_j.$ Then,
$\omega_j^{-\varepsilon(j)}.c_j \not \in {\rm acl}(\Omega_{j-1,\iota(j)}^{-\varepsilon(j)} \cup \{\omega_j^{-\varepsilon(j)}\})$ as required. Indeed, this follows from the decomposition
$${\rm acl}(\Omega_{j-1,\iota(j)}^{-\varepsilon(j)} \cup \{\omega_j^{-\varepsilon(j)}\}) = \left( {\rm acl}(\Omega_{j-1,\iota(j)}^{-\varepsilon(j)} \cup \{\omega\}) \setminus \{\omega\}\right) \cup \{\omega_j^{-\varepsilon(j)}\}.
$$
This completes the proof of (ii), since again in the last step there are infinitely many choices for $\beta_n=\omega_{l}^{-\varepsilon(l)}.c_l\notin{\rm acl}(\{\alpha_1,\beta_1,\ldots,\alpha_n\})$ as required.

Similarly, we need a special argument to finish the proof of (iii).
Indeed, for $j\in\{1,\ldots,l-1\}$, in Case~3 of the proof, we need that 
the set of all $x\in X$ such that $$x.c_j\notin{\rm acl}(\Omega_{j-1,\iota(j)}^{-\varepsilon(j)}\cup\{x\})$$ intersects with the infinite definable set $\Sigma_j$. By definition, this is the case if $c_j$ is not slender. The case $j=l$ is dealt with as above.

This completes the argument for (i), (ii) and (iii) in step $n\geq1$ of the induction. For the next step $n+1$, the proof is the same, but we replace $\Lambda_k^\varepsilon$ by the set $\Omega_k^\varepsilon$ from step $n$.
\end{proof}

We record the following three corollaries:

\begin{corollary}
A closed oligomorphic subgroup of ${\rm Sym}(\omega)$ does not admit strong mixed identities.
\end{corollary}

\begin{corollary}
Every mixed identity of the automorphism group of an $\omega$-categorical structure which is algebraically convex must have a small critical constant.    
\end{corollary}

\begin{corollary}
Every mixed identity of the automorphism group of an $\omega$-categorical structure must have a slender critical constant.    
\end{corollary}

\begin{corollary}\label{cor:small}
    If the automorphism group of an $\omega$-categorical structure has no non-trivial small elements,
    then it has no mixed identity.  
\end{corollary}

\begin{remark}
Despite numerous examples in Section \ref{sec:examples} we do not know whether the existence of non-trivial small or slender automorphisms for some $\omega$-categorical relational structure $(X,R)$ implies the existence of a mixed identity for $\Aut(X,R)$. 
\end{remark}

\section{Examples with mixed identities}
\label{sec:examples}
Now we will provide some examples of homogeneous relational structures with oligomorphic automorphism group for which we can apply Theorem~\ref{thm:main}. We will also provide mixed identities in the examples in this section; needless to say, those mixed identities are singular in accordance with Conjecture \ref{conj:autq}.

\subsection{The countably infinite set}\label{subsec:ctbl_st}

The simplest example of a homogeneous structure with an oligomorphic automorphism group is a \emph{countably infinite set} $X$ equipped with no relations at all. In this case, $(X,\varnothing)$ has no algebraicity, and its automorphisms are small if and only if they are slender if and only if they have finite support. Thus, we can apply Theorem~\ref{thm:main} to deduce that every mixed identity for ${\rm Sym}(X)$ has a critical constant with finite support, a fact which was originally proved by Hull--Osin in \cite{hullosin2016transitivity}.
Note that for $\tau$ a transposition, $w(x)=[\tau,x]^6$ is such an identity and there are many more of the same kind. Despite the abundance of mixed identities, we have the following consequence of Theorem~\ref{thm:main}:

\begin{theorem} \label{thm:nonnsmix}
There are no non-singular mixed identities for ${\rm Sym}(X)$.    
\end{theorem}

\begin{proof}
The proof proceeds as the proof of \cite{bradfordschneiderthom2023non}*{Corollary~1}. Let $w \in {\rm Sym}(X) \ast \mathbf{F}_r$ be a non-singular word with constants.
We replace the critical constants of finite support in $w$ one by one with the identity element. This leads to cancellation of the variables, changes of the neighboring constants, and reduces the word length in each step by at least two. Hence, this process ends after finitely many steps and leaves us with $w' \in {\rm Sym}(X) \ast \mathbf{F}_r$ without critical constants of finite support. It is clear that the content is unchanged in this process, i.e., $\epsilon(w)=\epsilon(w').$ In particular, the word length of $w'$ is positive. Moreover, it is clear from the construction that there exists a constant $f \in \mathbb{N}$ such that 
$$
{\rm dist}(w(h_1,\ldots,h_r),w'(h_1,\ldots,h_r))\leq f
$$
for all $h_1,\ldots,h_r \in {\rm Sym}(X),$ where ${\rm dist}$ denotes the Hamming distance of permutations.
Indeed, $f$ can be taken to be the sum of the support sizes of all critical constants that have been replaced. Now, by Theorem~\ref{thm:main}(ii), there exists a particular choice $h_1,\ldots,h_r \in H$ such that the support of $w'(h_1,\ldots,h_r)$ is larger than $f+1.$ This implies that $w(h_1,\ldots,h_r) \neq {\rm id}_X$ and the proof is complete.
\end{proof}

\subsection{The equivalence relation with $k$-element classes}\label{sec:eq}

Let $X$ be a countably infinite set and let $k \in {\mathbb N}_+$. On $X$ there is an \emph{equivalence relation $E_k\subseteq X^2$ with $k$-element equivalence classes}; note that 
this generalizes the example from the previous section, because $\Aut(X,E_1) = {\rm Sym}(X)$.
The structure $(X,E_k)$ is homogeneous and ${\rm acl}(Y)= \bigcup_{y \in Y} {\rm acl}(\{y\})$, but the structure is not algebraically convex. 
The automorphism group of $(X,E_k)$ is isomorphic to $S_k\wr{\rm Sym}(X/E_k)$, and the slender automorphisms are those in $S_k\wr{\rm Sym}_{\rm fin}(X/E_k)$, where ${\rm Sym}_{\rm fin}(X/E_k)$ denotes the subgroup of permutations of finite support. Theorem~\ref{thm:main}(iii) implies that every mixed identity for $S_k \wr {\rm Sym}(X/E_k)$ has a critical constant in $S_k \wr {\rm Sym}_{\rm fin}(X/E_k)$. While there are non-trivial singular mixed identities as $S_k^{X/E_k}$ is a normal subgroup of $\Aut(X,E_k)$ and this group has many such identities, we obtain the following as 
a direct consequence of Theorem~\ref{thm:nonnsmix}. 

\begin{corollary} The group $\Aut(X,E_k)$ does not have any non-singular mixed identites.
\end{corollary}

\subsection{The countably infinite dimensional vector space}
\label{sec:vectspace}
The \emph{countably infinite dimensional vector space} $\mathbb{F}_q^{\oplus \omega}$ over the finite field $\mathbb{F}_q$ is another interesting case. Further complications arise if we add a non-degenerate symmetric or alternating bilinear form $f \colon \mathbb{F}_q^{\oplus \omega} \times \mathbb{F}_q^{\oplus \omega} \to \mathbb{F}_q$ to the structure. We denote the respective automorphism groups by ${\rm GL}(\mathbb{F}_q^{\oplus \omega}), {\rm O}(\mathbb{F}_q^{\oplus \omega}),$ and ${\rm Sp}(\mathbb{F}_q^{\oplus \omega})$. {They are oligomorphic; see, e.g.,~\cite{evans1991small}.} 

A concrete relational structure for $\mathbb{F}_q^{\oplus \omega}$ such that the automorphism group is homogeneous is given by the infinite set of relations 
$$
R_{\lambda_1,\ldots,\lambda_n} \coloneqq \left\{(x_1,\ldots,x_n) \in (\mathbb{F}_q^{\oplus \omega})^n \,\middle\vert\, \sum_{i=1}^n \lambda_i x_i = 0 \right\}
$$
for $\lambda_1,\ldots,\lambda_n \in \mathbb{F}_q$ and $n \in \mathbb{N}.$
In the presence of a bilinear form $f$, we also add the relations $B_{\lambda} \coloneqq  \{(x,y) \in (\mathbb{F}_q^{\oplus \omega})^2 \mid f(x,y)=\lambda \}$ for $\lambda \in \mathbb{F}_q.$

The algebraic closure coincides with the linear hull for these examples. Note that these structures are not algebraically convex so that Theorem~\ref{thm:main}(ii) cannot be applied. Instead, we will apply Theorem~\ref{thm:main}(iii). Before we do so, let's characterize slender automorphisms in more familiar terms.

\begin{lemma}
An element $c\in {\rm GL}(\mathbb{F}_q^{\oplus \omega})$, ${\rm O}(\mathbb{F}_q^{\oplus \omega})$, and ${\rm Sp}(\mathbb{F}_q^{\oplus \omega})$ is slender if and only if ${\rm rk}(c-\lambda{\rm id})<\infty$ for some $\lambda\in\mathbb{F}_q$.
\end{lemma}

\begin{proof}
Let's first describe the infinite definable subsets in sufficient detail. In the case of no bilinear form, the infinite definable sets contain the complement of a subspace of finite dimension; in particular, they are co-finite.
In the other two cases, further constraints arise from the bilinear forms encoded in the relations $B_{\lambda}$ for $\lambda \in \mathbb{F}_q$. In any case, an infinite definable set contains an affine subspace of finite co-dimension.

Assume that ${\rm rk}(c-\lambda{\rm id})$ is finite and set $g_\lambda\coloneqq c-\lambda{\rm id}$. 
Then $W_\lambda\coloneqq{\rm im}(g_\lambda)$ is finite and for all $x\in \mathbb{F}_q^{\oplus \omega}$, we get
$x.(c-\lambda {\rm id}) \in W_{\lambda}$ and hence $x.c \in \langle W_{\lambda},x \rangle = {\rm acl}(W_\lambda \cup\{x\})$. Hence, $c$ is slender. 

Conversely, if ${\rm rk}(c-\lambda{\rm id})=\infty$ for all $\lambda\in\mathbb{F}_q$, then $Y.g_\lambda^{-1}\leq X$ has infinite codimension for all finite-dimensional $Y$. Now, an infinite definable set contains an affine subspace of finite codimension; hence, it cannot be contained in the union of $q$ subspaces of infinite codimension. Whence, we conclude that $c$ is not slender.
\end{proof}

We can now apply Theorem~\ref{thm:main}(iii).

\begin{theorem}
Every mixed identity for the groups ${\rm GL}(\mathbb{F}_q^{\oplus \omega}), {\rm O}(\mathbb{F}_q^{\oplus \omega})$, and $ {\rm Sp}(\mathbb{F}_q^{\oplus \omega})$ contains a critical constant $c$ such that ${\rm rk}(c-\lambda{\rm id})<\infty$ for some $\lambda\in\mathbb{F}_q$. 
\end{theorem}

A stronger result was already proven along with quantitative forms in finite dimensions in \cite{bradfordschneiderthom2023non} for ${\rm GL}(\mathbb{F}_q^{\oplus \omega})$ and will appear in a forthcoming paper \cite{schneiderthom2024simple} in the other cases.
While it is well-known that the groups ${\rm GL}(\mathbb{F}_q^{\oplus \omega})$, ${\rm O}(\mathbb{F}_q^{\oplus \omega})$, ${\rm Sp}(\mathbb{F}_q^{\oplus \omega})$ admit singular mixed identities, we obtain the following consequence:

\begin{corollary}
The groups ${\rm GL}(\mathbb{F}_q^{\oplus \omega}), {\rm O}(\mathbb{F}_q^{\oplus \omega})$, and ${\rm Sp}(\mathbb{F}_q^{\oplus \omega})$ do not admit non-singular mixed identities.
\end{corollary}

\begin{proof}
Note that Theorem~\ref{thm:main}(iii) implies that for every $n \in \mathbb{N}$ there exist $h_1,\ldots,h_{r}$ with $$\min\{{\rm rk}(w(h_1,\ldots,h_{r})-\lambda{\rm id}) \mid \lambda \in \mathbb{F}_q^{\times} \} \geq n.$$ 
The proof now proceeds as the proof Theorem~\ref{thm:nonnsmix} (see also \cite{bradfordschneiderthom2023non}*{Corollary~1}). 
\end{proof}

\begin{remark}
Note that ${\rm GL}(\mathbb{F}_q^{\oplus \omega})$ does not have non-trivial small elements, since any cofinite set of vectors spans the entire vector space. Hence, from this example we also see that Theorem~\ref{thm:main}(ii) does not hold without the assumption of algebraic convexity.
\end{remark}

\subsection{The countable unbounded dense linear order}\label{subsubsec:rat_ln}

Another interesting case to consider is the (up to isomorphism unique) \emph{countable unbounded dense linear order} $(\mathbb{Q},<)$, which is the Fra\"iss\'e limit of all  finite linear orders.

\begin{lemma}\label{lem:rats_no_alg_sml_aut}
The structure $(\mathbb{Q},<)$ has no algebraicity and its small automorphisms are precisely those that fix a non-trivial interval pointwise. 
\end{lemma} 



This follows from the fact that every infinite definable set contains a non-trivial interval. Then Lemma~\ref{lem:rats_no_alg_sml_aut} and Theorem~\ref{thm:main}(ii) implies the following. 

\begin{corollary}
Every mixed identity for $\Aut(\mathbb{Q},<)$ has a critical constant which fixes a non-trivial interval pointwise.
\end{corollary}

A variation of the proof of \cite{zarzycki2010limits}*{Proposition~2.2} shows that the group $\Aut(\mathbb{Q},<)$ has singular mixed identities.
We denote by ${\rm supp}(g)= \{x \in \mathbb Q \mid x.g \neq x\}$ the support of an element $g \in \Aut(\mathbb{Q},<)$.

\begin{lemma}
There exist singular mixed identities for $\Aut(\mathbb{Q},<)$.
\end{lemma}

\begin{proof}
Let $g_1, g_2, g_3 \in \Aut(\mathbb{Q},<)$ be three non-trivial elements having their support in non-trivial intervals $I_1=(a_1,b_1),I_2=(a_2,b_2),I_3=(a_3,b_3)$ with the property that $a_1<b_1<a_2<b_2<a_3<b_3$. Then,
$$
w\coloneqq [[g_1^x,g_3],[g_2^x,g_2]]\in\Aut(\mathbb{Q},<)\ast\langle x\rangle
$$
is a mixed identity for $\Aut(\mathbb{Q},<).$ To see this, we need two basic facts:
\begin{enumerate}[$(i)$]
\item ${\rm supp}(h^g)={\rm supp}(h).g$ for any $g,h \in \Aut(\mathbb{Q},<)$, and
\item $[g,h]={\rm id}_\mathbb{Q}$ if ${\rm supp}(g)\cap\sup(h)=\varnothing.$
\end{enumerate}

If $I_1.g\cap I_3=\varnothing$, then by (i) ${\rm supp}(g_1^g)\subseteq I_1.g$ and by (ii) we conclude $[g_1^g,g_3]={\rm id}_\mathbb{Q}$ and hence $w(g)={\rm id}_\mathbb{Q}.$ On the other side, if $I_1.g\cap I_3\neq\varnothing$ then by the monotonicity of $g$, we must have $I_2.g\cap I_2=\varnothing$. Hence, we conclude in a similar way $[g_2^g,g_2]={\rm id}_\mathbb{Q}$ and hence $w(g)={\rm id}_\mathbb{Q}.$ This finishes the proof.
\end{proof}

We will come back to the study of non-singular mixed identities for the group $\Aut(\mathbb{Q},<)$ in Section \ref{sec:autQ}.

\subsection{The generic cyclic order}

Consider now $\Aut(\mathbb Q/\mathbb Z,t)$, where $t$ is the ternary betweenness relation, i.e., $t(\alpha,\beta,\gamma)$ holds if and only if $\beta$ lies in the counter-clockwise oriented segment from $\alpha$ to $\gamma$, where we view $\mathbb{Q}/\mathbb{Z}$ as a subset of ${\mathbb S}^1$ in a natural way. We can also consider the group $\Aut^{\pm}(\mathbb Q/\mathbb Z,t)$ of automorphisms of $\mathbb Q/\mathbb Z$ that either preserve or reverse the betweenness relation. This group contains $\Aut(\mathbb Q/\mathbb Z,t)$ as a subgroup of index $2$ and acts 3-transitively on $\mathbb Q/\mathbb Z$. The result follows also from work of Le Boudec--Matte Bon; see \cite{leboudecmattebon2022triple}*{Lemma~4.6}.

\begin{theorem}\label{thm:Q}
Let $g_1,g_2,g_3,g_4,g_5 \in \Aut(\mathbb Q/\mathbb Z,t)$ be five non-trivial elements having their support in cyclically counter-clockwise ordered, pairwise disjoint intervals
$I_1,I_2$, $I_3,I_4,I_5 \subseteq\mathbb Q/\mathbb Z.$ Then,
$$w(x) \coloneqq [[[g_1^x,g_4],[g_3^x,g_5]],[g_2^x,g_2]]$$
is a mixed identity for $\Aut(\mathbb Q/\mathbb Z,t)$ and $$w'(x)\coloneqq w(x^2)$$ 
is a mixed identity for $\Aut^{\pm}(\mathbb Q/\mathbb Z,t).$
\end{theorem}
\begin{proof}
The proof uses the same idea as above and the geometric fact that for $g \in \Aut(\mathbb Q/\mathbb Z,t)$ the condition
$I_1.g \cap I_4 \neq \varnothing$ together with $I_3.g \cap I_5 \neq \varnothing$ implies that $I_2.g \cap I_2 = \varnothing.$ Thus, $[g_2^g,g_2]$ is trivial unless $[g_1^g,g_4]$ or $[g_3^g,g_5]$ are trivial. We conclude that $w(g)={\rm id}_{\mathbb Q/\mathbb Z}$ for all $g \in \Aut(\mathbb Q/\mathbb Z,t).$ Thus, $w(x)$ is a mixed identity for $\Aut(\mathbb Q/\mathbb Z,t)$ and, clearly, $w'(x) = w(x^2)$ is a mixed identity for $\Aut^{\pm}(\mathbb Q/\mathbb Z,t)$.
\end{proof}

\subsection{The generic semi-linear order with joins}\label{sec:semilinear}

{A \emph{semilinear order}~\cite{macpherson2011survey} is a partially ordered set 
$(S,<)$ such that
\begin{itemize}
    \item for each $a \in S$
the set $\{x \in S \mid a<x\}$ is linearly ordered by $<$, and
\item for all $x,y \in S$ there exists $z \in S$ such that $x<z$ and $y<z$. 
\end{itemize}
}
Countable semilinear orders where all maximal chains are isomorphic to $({\mathbb Q},<)$ have been classified by Droste~\cite{droste1985structure}. 
There exists such a semilinear order $(S,<)$ where for any two $x,y$ 
\begin{itemize}
\item there exists a supremum of $x$ and $y$ in $S$, i.e., a (unique) smallest upper bound $\sup(x,y)$ for $x$ and $y$, and 
\item for every $x \in S$, the equivalence relation $E_x$ on $T_x \coloneqq\{z \in S \mid z \leq x \}$ which contains all pairs $(u,v) \in T_x^2$ such that there exists $w \in T_x$ with $u,v < w < x$, has infinitely many equivalence classes.   
\end{itemize}
The semilinear order $(S,<)$ is unique up to isomorphism, and has an oligomorphic automorphism group, but is not homogeneous. However, as explained in Remark~\ref{rem:homogen} it has a homogeneous expansion with the same automorphism group. In model theory, it is a frequently studied example in the context of structures that are NIP but not stable; for recent results, see~\cite{kaplanrzepeckisiniora2021automorphism}. 
In the following, we call $(S,<)$ the \emph{generic semilinear order with joins}.  

\begin{lemma}\label{lem:smi_lin_ord}
The generic semi-linear order with joins $(S,<)$ does have algebraicity and is algebraically convex. The small elements $c\in\Aut(S,<)$ are those fixing a non-trivial interval.
\end{lemma}
 
\begin{proof}
For the first two claims, let $Y$ be a finite subset of $S$. Note that the algebraic closure of $Y$ is the smallest subset of $S$ containing $Y$ which is closed under joins. Moreover, it is easy to see that 
$\acl(Y \cup \{x\}) = \acl(Y) \cup \{\bar x,x\}$ for some $\bar x \in S$ (we may have $x = \bar x$). 
Note that  
 $\acl(Y) \cup \{\bar x\}$ is algebraically closed. In particular, $(S,<)$ is algebraically convex.


For the third claim, one can prove that each infinite  definable set contains a non-trivial interval $(a,b)=\{x\in S\mid a<x<b\}$ ($a<b$, $a,b\in S$). Hence $c\in\Aut(S,<)$ is small if and only if it fixes such an interval pointwise (in analogy to Section~\ref{subsubsec:rat_ln} above).
\end{proof}

We obtain from Lemma~\ref{lem:smi_lin_ord} and Theorem~\ref{thm:main}:

\begin{corollary}
Every mixed identity of $\Aut(S,<)$ contains a critical constant fixing a non-trivial interval.
\end{corollary}

We can also prove the existence of mixed identities. The following construction is inspired by a result in~\cite{leboudecmattebon2022triple}.

\begin{lemma}
There exist mixed identities for $(S,<)$.
\end{lemma}

\begin{proof}
Suppose that $g,h\in\Aut(S,<)$ are supported in disjoint half-trees (these are sets of the form $\{x\in S\mid x<y\}$). In order to see that such $g$ and $h$ exist, we consider two incomparable elements $y,z\in S$. Let $y_1$ and $y_2$ be incomparable below $y$. Then, there exists an automorphism $g'$ that fixes $y$ and interchanges $y_1$ and $y_2$. Define $g\in\Aut(S,<)$ to fix $y$, everything above $y$ and everything incomparable to $y$; and to be equal to $g'$ below $y$. This defines an automorphism of the generic semi-linear order. Similarly, you find $h\in\Aut(S,<)$ which only acts below $z$. Now, $g$ and $h$ are supported on the disjoint half-trees $\{x\in S\mid x<y\}$ and $\{x\in S\mid x<z\}$. Consider an arbitrary element $k\in\Aut(S,<)$. There are only three possible cases. In each of them we use that two automorphisms with disjoint support commute:
\begin{enumerate}[$(i)$]
\item If $y.k$ is incomparable to $y$, then $g^k$ and $g$ commute.
\item If $y.k$ is below $y$, then $z$ is incomparable to $y.k$, and $g^k$ and $h$ commute. 
\item If $y.k$ is greater than or equal to $y$, then $z.k$ is incomparable to $y$, and so $h^k$ and $g$ commute.
\end{enumerate}
We conclude that 
$$
w(x)=[[g^x,g],[[h^x,g],[g^x,h]]]\in\Aut(S,<)\ast\langle x\rangle
$$ 
is a mixed identity for the automorphism group.
\end{proof}

\subsection{The Wa\.{z}ewski dendrite}\label{sec:dendrite}

A \emph{dendrite} $X$ is a compact, metrizable, connected and locally connected topological space without simple closed curves. Then, for each pair of distinct points $x,y \in X$, there exists a unique simple path from $x$ to $y$ in $X.$ We say that a set $W \subseteq X$ is \emph{arcwise dense}, if it intersects with each such simple path.
See \cite{duchesnemonod2019structural} for more background and references. 

For a dendrite $X$, we call $x \in X$ a \emph{branch point}, if $X \setminus \{x\}$ consists of at least three connected components. The cardinality of the set of connected components of $X \setminus \{x\}$ is called the \emph{order} of $x.$ We denote by ${\rm Br}(X)$ the set of branch points of $X$.

Throughout, let $S \subseteq \{3,4,\ldots\} \cup \{\infty\}$ be a finite subset. The Wa\.{z}ewski dendrite $D_S$ is uniquely characterized by the property that each branch point has order in $S$ and that the set of branch points of order $s \in S$ is arcwise dense in $D_S$. It follows that the set ${\rm Br}(D_S)$ is countable.

It was proved in \cite{duchesnemonod2019structural}*{Proposition 2.4} that the natural group homomorphism $$\varphi \colon {\rm Homeo}(D_S) \to\Sym({\rm Br}(D_S))$$ is a homeomorphism onto its image and it follows from \cite{duchesnemonod2019structural}*{Corollary~6.8} that its image is an oligomorphic permutation group.
A natural relation that is preserved by this permutation action is the relation $R$ on triples with $(x,y,z) \in R$ if and only if $y$ lies on the simple path from $x$ to $z$. Moreover, one can identify the image of $\varphi$ with $\Aut({\rm Br}(D_S),R).$
From this relation, one can define uniquely the midpoint of the tripod spanned by three points $x,y,z \in {\rm Br}(D_S)$. We see from this that the algebraic closure of three points can be non-trivial. More generally, the algebraic closure of a finite set $F \subseteq {\rm Br}(D_S)$ consists of the vertices of the finite tree spanned by $F.$

\begin{lemma}
The homogeneous relational structure $({\rm Br}(D_S),R)$ is algebraically convex. The small automorphisms are those fixing a non-trivial simple path point-wise.
\end{lemma}

\begin{corollary}
Every non-trivial mixed identity of $\Aut({\rm Br}(D_S),R)$ contains a critical constant that fixes a non-trivial simple path.
\end{corollary}

We say that $A \subseteq {\rm Br}(D_S)$ is a \emph{half-tree} if there exists a point $x \in D_S \setminus {\rm Br}(D_S)$ of order $2$ and $A = {\rm Br}(D_S) \cap B$, where $B$ is one of the two connected components of $D_S \setminus \{x\}.$ Using \cite{leboudecmattebon2022triple}*{Lemma~3.5}, we can now construct a mixed identity for the group $\Aut({\rm Br}(D_S),R)={\rm Homeo}(D_S).$ The following is restating this lemma in our setting.

\begin{lemma}
Let $A_1,A_2,A_3 \subseteq {\rm Br}(D_S)$ be disjoint half-trees and $g \in \Aut({\rm Br}(D_S),R).$ Then either
\begin{enumerate}
\item there exists $i \in \{1,2,3\}$ with $A_i$ and $A_i.g$ disjoint, or
\item $A_i.g=A_i$ for all $i \in \{1,2,3\}.$
\end{enumerate}
\end{lemma}

We can now follow the arguments in \cite{leboudecmattebon2022triple}*{Lemma~3.6} and the discussion after its proof to obtain the following corollary:

\begin{corollary}
Let $h_1,h_2,h_3 \in \Aut({\rm Br}(D_S),R)$ be automorphisms that have support in disjoint half-trees and let $g \in \Aut({\rm Br}(D_S),R)$ be arbitrary. Then, either $[h_i^g,h_i]=1$ for one $i \in \{1,2,3\}$ or $[h_i^g,h_j] = 1$ for all $i,j \in \{1,2,3\}$ with $i \neq j.$ In particular, the word
$$w(x) = [[[h_1^x,h_1],[h_2^x,h_2]],[[h_3^x,h_3],[h_1^x,h_2]]]$$
is a mixed identity for $\Aut({\rm Br}(D_S),R)$.
\end{corollary}

\section{Examples without mixed identities}
This section contains examples of oligomorphic automorphism groups for which there do not exist any mixed identities at all.

\subsection{The Rado graph}\label{sect:rado}
The celebrated \emph{Rado graph} $R=(V,E)$ is the Fra\"iss\'e limit of the class of all finite graphs. Here $V$ is the set of  vertices and $E$ the set of edges. It is well  known and easy to see that the Rado graph is the up to isomorphism unique countable graph with the so-called \emph{extension property}, i.e., for all finite disjoint vertex subsets $A,B \subseteq V$ there exists $x \in V$ such that $x$ is adjacent to all vertices in $A$ and to no vertex in $B$. Even though the following lemma is subsumed by a more general result by Macpherson and Tent (Theorem~\ref{thm:macphersonTent}) we also present a proof because it is instructive and simple.

\begin{lemma}\label{lem:un_rado}
The Rado graph $R = (V,E)$ has no algebraicity and no non-trivial small automorphisms.
\end{lemma}

\begin{proof}
The class of all finite graphs even has the free amalgamation property, and hence the strong amalgamation 
property, and so has no algebraicity by Remark~\ref{rem:sap}.
To establish that $R$ has no non-trivial small automorphisms, we prove Property~(ii) of Lemma~\ref{lem:no_smll_elts}. Let $Y \subseteq V$ be finite and $x\in V \setminus Y$. Set $Z\coloneqq x.\stab(Y)$. Consider $c\in\stab(Z).$ Let $u,v\in V\setminus Z$ be arbitrary distinct vertices. At first we show that $\supp(c)\subseteq Y\cup Y.c^{-1}$. Assume that $u,v\notin Y$. By the extension property of the Rado graph applied to $Y\cup\{u,v\}$, there exists a vertex $z\in V$ that behaves to $Y$ as $x$ does, and such that $u$ is connected to $z$, but $v$ is not. Thus by homogeneity, we have that $z\in x.\stab(Y)$ and hence $z.c=z$. But then $u$ cannot be mapped to $v$ by $c$ by the choice of $z$. 

Hence, we conclude that $c$ has finite support. Let $w\in\supp(c)$. By the defining property, there exists $y\in V\setminus\supp(c)$ which is adjacent to $w$ but to none of the vertices in $\supp(c)\setminus\{w\}$. Then the edge $(y,w)$ is mapped to the non-edge $(y.c,w.c)=(y,w.c)$ and we obtain a contradiction. Hence, $c$ must be trivial.
\end{proof}

Lemma~\ref{lem:un_rado} and 
Corollary~\ref{cor:small}
imply:

\begin{corollary}\label{cor:rado}
The group $\Aut(R)$ has no mixed identities.
\end{corollary}



\subsection{Free amalgamation}
Examples of homogeneous structures whose age has the free amalgamation property are the Rado graph from the previous section, the random $k$-uniform hypergraph, for some $k \geq 3$, the homogeneous graph whose age consists of all finite triangle-free graphs (or, more generally, $K_n$-free graphs, for some fixed $n \geq 3$), the uncountable family of homogeneous directed graphs constructed by Henson~\cite{henson1972countable}, and many more.
Already for the random triangle-free graph, the easy proof that we presented for the Rado graph in the previous section fails, since for some choices of $Y$ and $x$, there might not exist a vertex $z \in V$ as it does not contain triangles. 
All of the listed examples above also satisfy the additional assumptions made in the following theorem. 

\begin{theorem}[Corollary 2.10 in \cite{MacphersonTentSimplicity}]
\label{thm:macphersonTent}
    Let $(X,R)$ be a homogeneous relational structure whose age has the free amalgamation property such that
$\Aut(X,R)$ is transitive and 
$\Aut(X,R) \neq {\rm Sym}(X)$. 
Then  $(X,R)$ has no non-trivial small automorphisms. 
\end{theorem}

Combined with Corollary~\ref{cor:small} we now obtain the following.

\begin{theorem}\label{thm:free}
Let $(X,R)$ be a homogeneous relational structure whose age has the free amalgamation property such that
$\Aut(X,R)$ is transitive and 
$\Aut(X,R) \neq {\rm Sym}(X)$. 
Then $\Aut(X,R)$ has no mixed identities.
\end{theorem}

\subsection{The random poset} 

The \emph{random poset} $(P,<)$ is defined as the Fra\"iss\'e limit of the class of all finite posets; note that this class does not have the free amalgamation property, but the strong amalgamation property, and hence $(P,<)$ has no algebraicity. For this relational structure the following holds.

\begin{lemma}\label{lem:un_pset}
The random poset $(P,<)$ has no non-trivial small automorphisms.
\end{lemma}

\begin{proof}
Let $a,b \in P$ with $a<b$. We will prove that every automorphism that fixes $Z\coloneqq\{x\in P\mid a<x<b\}$ is trivial. We will write $x\perp y$ if neither $x<y$ nor $y<x,$ i.e., $x$ and $y$ are incomparable.

Suppose that $c$ and $c'$ are incomparable to $b$ and both larger than $a$. Then, by homogeneity, there exists $x \in P$ with $a<x<b$, $x<c$ but $x \not < c'.$ Now, any automorphism that fixes $Z$ cannot map $c$ to $c'$. It follows that $\stab(Z)$ fixes $\{c \mid c \perp b, a<c\}$ pointwise.
 
Similarly, $\stab(Z)$ fixes $\{c \mid c < b, a \perp c\}$ and $\{c \mid c \perp b, c\perp a\}$ pointwise. Consider now $c$ and $c'$ with $b<c$ and $b<c'$. By homogeneity, there exists $x \in P$ incomparable to $b$ such that $x<c$ but $x \not < c'.$ Since $x$ is fixed by $\stab(Z)$ by the previous arguments, no element of $\stab(Z)$ can map $c$ to $c'$. It follows that $\stab(Z)$ fixes the set $\{c \mid b<c\}$ pointwise. Similarly, it fixes $\{c \mid c<a\}$ pointwise. This covers all cases and implies that $\stab(Z)$ is trivial. Hence the claim is proven.
\end{proof}

As a consequence of Lemma~\ref{lem:un_pset} 
and Corollary~\ref{cor:small}
we obtain:

\begin{corollary}
The group $\Aut(P,<)$ has no mixed identities.
\end{corollary}

\subsection{The random permutation}
The \emph{random permutation} is the (up to isomorphism unique) countable homogeneous structure such that every finite structure with two linear orders is isomorphic to a substructure of $(P,<_1,<_2)$. Note that $(P,<_1)$ and $(P,<_2)$ are each isomorphic to $({\mathbb Q},<)$, so we may identify $P$ with a subset of $\mathbb{Q} \times \mathbb{Q}$, together with $<_1$ and $<_2$ induced by the two orderings coming from the coordinate projections.

\begin{remark} 
Just like $({\mathbb Q};<)$, the first-order theory of the random permutation $(P;<_1,<_2)$ is known to be NIP (see, e.g.~\cite{SimonRank1}).
However, its automorphism group behaves differently from $\Aut({\mathbb Q};<)$ with respect to the existence of mixed identities, as we will see below.
\end{remark}

\begin{lemma}
    $\Aut(P,<_1,<_2)$ does not have small automorphisms.
\end{lemma}
\begin{proof} Let $Y\subseteq P$ be finite, let $x \in P \setminus Y$, 
and let $Z\coloneqq x.\stab(Y)$; clearly, $Z$ is infinite and definable in $(P;<_1,<_2)$. 
Consider $\alpha \in \stab(Z)$. 
It suffices to show that $\alpha = \text{id}_X$, because every infinite definable subsets of $(P;<_1,<_2)$ must contain a set of the form $Z$ as described above (see the proof of Lemma~\ref{lem:no_smll_elts} implication (ii)$\Rightarrow$(iii)). Let $a,b \in P$ be distinct. 
We may assume that $a <_1 b$.
First we consider the case that 
there exists $c \in Z$ such that $a<_1 c <_1 b$. In this case, $a.\alpha <_1 c.\alpha = c <_1 b$ and we thus have that $a.\alpha \neq b$. Therefore, 
$\alpha|_{Z_1} = \id_{Z_1}$ where
$Z_1$ is the set of all elements that lie in the same orbit as $x$ with respect to $\Aut(P;<_1)$. Likewise, $\alpha|_{Z_2} = \id_{Z_1}$ 
where $Z_2$ is the set of all elements that lie in the same orbit as $x$ with respect to $\Aut(P;<_2)$. In the general case, there exists $c \in Z_2$ such that $a <_1 c <_1 b$, 
and we conclude that $\alpha = \id_P$, as required. 
\end{proof} 

\begin{corollary}
$\Aut(P,<_1,<_2)$ has no mixed identities.
\end{corollary}

\subsection{Parametrised families of linear orders}
Consider $k \in \mathbb{N}$ and consider the class $\mathcal R_k$ of finite structures $(X,R)$, where $X$ is a finite set and $R$ is a $(k+2)$-ary relation with the properties
\begin{enumerate}
\item $(x_1,\ldots,x_k,a,b) \in R$ only if $x_1,\ldots,x_k,a,b$ are pairwise distinct,
\item for a fixed $k$-tuple $x=(x_1,\ldots,x_k) \in X^k$ of pairwise distinct elements, the binary relation
$$L_k(x) \coloneqq \;  \{(a,b) \in (X \setminus \{x_1,\ldots,x_k\})^2 \mid (x_1,\ldots,x_k,a,b) \in R \}$$ 
defines a linear order on $X(x) \coloneqq X \setminus \{x_1,\ldots,x_k\}.$
\end{enumerate}

\begin{lemma}
The class $\mathcal R_k$ is closed under substructures, isomorphisms, and has the strong amalgamation property.
\end{lemma}

\begin{proof}
Only the strong amalgamation property requires a proof. Let $\bA, \bB_1, \bB_2 \in {\mathcal R}_k$ be such that $A = B_1 \cap B_2$ and $\bA$ is a substructure of both $\bB_1$ and $\bB_2$. We need to find a structure $\bC \in {\mathcal R}_k$ and isomorphisms $e_i$, for $i \in \{1,2\}$, 
between $\bB_i$ and a substructure of $\bC$ such that $e_1(a) = e_2(a)$ for every $a \in A$. 

The domain of $\bC$ will be $C \coloneqq B_1 \cup B_2$. 
For every $x \in C^k$,
we define a linear order $L$ 
on $C \setminus \{x_1,\dots,x_{k}\}$ as follows. 
If $x \in (B_i)^k$, for $i \in \{1,2\}$, then 
$L_i$ denotes the restriction of the linear order $L_k(x)$ to $B_i \setminus \{x_1,\dots,x_{k}\} $. 
If $x \in A^k$, then 
$L$ is any linear extension of $L_1 \cup L_2$. 
Otherwise, if $x \in (B_i)^k$
for $i \in \{1,2\}$, then $L$ is any linear extension of $L_i$ to all of $C \setminus \{x_1,\dots,x_{k}\}$. 
Finally, if $x$ is neither in $(B_1)^k$  nor in $(B_2)^k$, 
then $L$ is any linear order on $C \setminus \{x_1,\dots,x_{k}\}$. 
The relation $R$ of $\bC$ 
then consists of all tuples of the form $(x,a,b)$ for $(a,b) \in L$; clearly, $\bC \in {\mathcal R}_k$. It is also clear that the identity map on $B_i$ is an isomorphism between $\bB_i$ and a substructure of $\bC$, which proves the strong amalgamation property. 
\end{proof}

Consider now the Fra\"{\i}ss\'{e} limit $\bX_k=(X_k,R_k)$ of the class $\mathcal R_k$ and its automorphism group $G_k=\Aut(\bX_k)$ (see, e.g.,~\cite{hodges1997shorter}). Note that the action of $G_k$ on $X_k$ is $(k+1)$-transitive and oligomorphic.

Let $x=(x_1,\ldots,x_k)$ be a $k$-tuple of distinct elements of $X_k$. 
The homogeneity of $(X_k,R_k)$ implies that $(X_k(x),L_k(x))$ is dense and unbounded, and hence isomorphic to $(\mathbb Q,<)$. We write $y<_x z$ for $(y,z) \in L_k(x).$

\begin{lemma}
Every infinite  definable set contains a set of the form $\{y \in X_k \mid a<_x y <_x b \}$ for some $k$-tuple $x=(x_1,\ldots,x_k)$ of distinct elements of $X_k$ and some $a,b \in X_k(x)$ with $a<_x b$.
\end{lemma}

\begin{lemma} 
The structure $\bX_k$ has no algebraicity, and for $k \geq 1$, it has no non-trivial small automorphisms.
\end{lemma}

\begin{proof}
No algebraicity follows from the strong amalgamation property and Remark \ref{rem:sap}. It remains to exclude the existence of small automorphisms. 
Let $Z$ be an infinite set which is definable in $\bX_k$
and let $c \in \stab(Z)$. 
Note that $Z$  contains a set of the form
$$
\bigcap_{j=1}^d \left\{y \in X_k \mid a_j<_{x^{(j)}} y <_{x^{(j)}} b_j \right\}
$$
for some $d \in {\mathbb N}$, some $k$-tuples $x^{(j)}=(x^{(j)}_1,\ldots,x^{(j)}_k)$ of distinct elements of $X_k$, and some $a_j,b_j \in X_k(x^{(j)})$ with $a_j<_{x^{(j)}} b_j$.

Let now $z_1,z_2 \in \bigcap_{j=1}^d X_k(x^{(j)}).$ By homogeneity and since $k \geq 1$, there exists a tuple $x' = (x'_1,\ldots,x'_k) \in (X_k)^k$  satisfying
\begin{itemize}
    \item $a_j<_{x^{(j)}} x'_i <_{x^{(j)}} b_j$ for all $i \in \{1,\ldots,k\}$ and $j \in \{1,\ldots,d\}$, \item 
    $z_1 <_{x'} a$,  and
    \item $ a <_{x'} z_2$. 
\end{itemize}
    Note that $x' \in Z^k$. Thus $x'.c=x'$, 
    and since 
    $z_1 <_{x'} z_2$ we have 
    $z_1.c \neq z_2$. This implies that $c$ fixes every element in $\bigcap_{j=1}^d X_k(x^{(j)})$. Let now $y =(y_1,\ldots,y_k)$ be a $k$-tuple of distinct elements of $\bigcap_{j=1}^d X_k(x^{(j)})$ satisfying $a_j<_{x^{(j)}} <y_i<_{x^{(j)}} b_j$. 
    Then $y.c=y$ and the linear order 
    $L_k(y)$ 
    on $\bigcup_{j=1}^d \{x^{(j)}_1, \ldots, x^{(j)}_k\}$ needs to be preserved by $c$. Hence, $x^{(j)}_i.c=x^{(j)}_i$ for all $i \in \{1,\dots,k\}$ and $j \in \{1,\dots,d\}$. This proves that $c$ is trivial and completes the proof.
\end{proof}

\begin{corollary}
For $k \geq 1$, 
the group $\Aut(\bX_k)$ has no mixed identities.
\end{corollary}

\section{Ruling out non-singular mixed identities for ${\rm Aut}(\mathbb{Q},<)$} \label{sec:autQ}

In this section, we continue to consider the group $\Aut(\mathbb{Q},<)$ of automorphisms of $(\mathbb{Q},<)$, see also Section \ref{subsubsec:rat_ln}. Our aim is to prove Conjecture~\ref{conj:autq} under some assumption on the content of the word with constants.

\subsection{${\rm Homeo}^+([0,1])$}

Our technique applies to the slightly larger group ${\rm Homeo}^+([0,1])$ of orientation-preserving homeomorphisms of the closed interval $[0,1] \subseteq {\mathbb R}$, equipped with the metric topology. We will need the following folklore result.

\begin{lemma}\label{lem:dense} 
There exists an injective continuous homomorphism $$\varphi\colon\Aut(\mathbb{Q},<) \to {\rm Homeo}^+([0,1])$$ with dense image.
\end{lemma}

Every mixed identity of $\Aut(\mathbb{Q},<)$ extends by continuity to a mixed identity of ${\rm Homeo}^+([0,1])$ and every mixed identity of ${\rm Homeo}^+([0,1])$ with constants in $\Aut(\mathbb{Q},<)$ restricts to a mixed identity of $\Aut(\mathbb{Q},<).$ We will therefore restrict attention to ${\rm Homeo}^+([0,1])$ from here on.

\vspace{0.2cm}

The main result of this section is the following theorem.

\begin{theorem}\label{thm:mixidaut}
Let $w \in{\rm Homeo}^+([0,1]) \ast \mathbf{F}_2$ be a word with constants such that $\epsilon(w) \not \in [[\mathbf{F}_2,\mathbf{F}_2],[\mathbf{F}_2,\mathbf{F}_2]]$. Then $w$ is not a mixed identity for ${\rm Homeo}^+([0,1])$.
\end{theorem}

We also consider the subgroup ${\rm Homeo}^+_{\rm Lip}([0,1])$ of bi-Lipschitz homeomorphisms, which serves as an important tool in the proof. The following result can be found as Th\'{e}or\`{e}me D in \cite{deroinkleptsynnavas2007dynamique}.

\begin{theorem}[Deroin--Kleptsyn--Navas]\label{thm:conjugate}
Every countable subgroup of ${\rm Homeo}^+([0,1])$ is conjugate to a subgroup of ${\rm Homeo}^+_{\rm Lip}([0,1])$.
\end{theorem}

Note that there is a natural action of the group ${\rm Homeo}^+([0,1])$ on $L^{\infty}([0,1],\mathbb R)$ by composition. This turns $L^{\infty}([0,1],\mathbb R)$ into a left-module over the integral group ring of the group ${\rm Homeo}^+([0,1])$ by the formula $(g\xi)(t) = \xi(t.g),$ for $g \in {\rm Homeo}^+([0,1])$ and $\xi \in L^{\infty}([0,1],\mathbb R).$

Recall that a cocycle $\mu \colon {\rm Homeo}^+([0,1]) \to L^{\infty}([0,1],\mathbb R)$ is a map that satisfies the equation
$$ 
\mu(gh) = g\mu(h) + \mu(g), \quad \forall g,h \in {\rm Homeo}^+([0,1]).
$$ An important example of a cocycle is the \emph{derivative cocycle} or \emph{Rademacher cocycle.} 

\begin{lemma}[Rademacher]
There exists a cocycle $$\mu \colon {\rm Homeo}^+_{\rm Lip}([0,1]) \to L^{\infty}([0,1],\mathbb R),$$ which is uniquely determined by the equation
$$
g(t) = \int_{0}^t \exp(\mu(g)(s))ds
$$
for all $g \in {\rm Homeo}^+_{\rm Lip}([0,1]).$
\end{lemma}

Note that $\mu(f) = \log\left(\frac{df}{dx} \right)$ if $f$ is continuously differentiable. Moreover, the same formula applies, when we note that a Lipschitz function is almost everywhere differentiable and interpret the terms in a suitable form. The cocycle equation $\mu(gh) = g \mu(h)+ \mu(g)$ is then a simple consequence of the chain rule for differentiation.

As a warm-up, we prove that there are no non-singular mixed identities in one variable. This already covers the case where  $\epsilon(w) \notin [\mathbf{F}_2,\mathbf{F}_2]$ by a simple reduction.

\begin{theorem} \label{thm:onevar}
The group ${\rm Homeo}^+([0,1])$ does not satisfy a non-singular mixed identity in one variable.
\end{theorem}

\begin{proof}
Suppose that 
$$w(x)=x^{\varepsilon(1)}h_1x^{\varepsilon(2)}\cdots x^{\varepsilon(l)}h_l$$
is a non-singular identity with constants in ${\rm Homeo}^+([0,1])$. We are assuming that $\varepsilon(i) \in \{\pm 1\}$ for all $i \in \{1,\ldots,l\}$.
Because $w$ is non-singular, we know that $e \coloneqq  \sum_{i=1}^{l} \varepsilon(i) \neq 0$ and we may assume without loss of generality that $e$ is positive.

Consider a countable subgroup $\Lambda$ of ${\rm Homeo}^+([0,1])$ containing $h_1,\ldots,h_l$. By Theorem~\ref{thm:conjugate}, we may assume without loss generality that $\Lambda$ is contained in ${\rm Homeo}^+_{\rm Lip}([0,1])$. 
We set
$C \coloneqq \max_{1 \leq i \leq l} \|\mu(h_i)\|_{\infty}$
so that $$k \mu(h_i) \geq -\|\mu(h_i)\|_{\infty} \geq - C$$ for every $k \in {\rm Homeo}^+([0,1])$. 
From the cocycle identity, we obtain that for every 
$g \in {\rm Homeo}^+_{\rm Lip}([0,1])$: 
\begin{eqnarray*}
 0 &=& \mu(w(g)) \\ &=& \sum_{i=1}^l g^{\varepsilon(1)}h_1 \cdots g^{\varepsilon(i-1)}h_{i-1} ( \mu(g^{\varepsilon(i)}) + g^{\varepsilon(i)} \mu(h_i)) \\
& \geq& -lC + \sum_{i=1}^l g^{\varepsilon(1)}h_1 \cdots g^{\varepsilon(i-1)}h_{i-1} \mu(g^{\varepsilon(i)}) \\
&=& -lC + \alpha \mu(g)
\end{eqnarray*}
with (using that $\mu(g^{-1}) = -g^{-1} \mu(g)$) $$\alpha \coloneqq  \sum_{i:e(i)=+1} g^{\varepsilon(1)}h_1 \cdots g^{\varepsilon(i-1)}h_{i-1} - \sum_{i:\varepsilon(i)=-1} g^{\varepsilon(1)}h_1 \cdots g^{\varepsilon(i-1)}h_{i-1} g^{-1} $$
considered as an element in the integral group ring of ${\rm Homeo}^+_{\rm Lip}([0,1])$. Note that the augmentation of $\alpha$ is equal to $e$.

For $\lambda\geq 1$, we consider now the specific bi-Lipschitz homeomorphism 
$$t.g_{\lambda} \coloneqq  \begin{cases} \lambda t & 0 \leq t \leq 1/(2\lambda) \\
\frac{1+ t-1/\lambda}{2-1/\lambda} & 1/(2\lambda) \leq t \leq 1. \end{cases}$$
 Then, the germ of the function $\alpha \mu(g_{\lambda})$ near $0$ is equal to $e \log(\lambda)$. But this is incompatible with the inequality $\alpha \mu(g_{\lambda}) \leq lC$ for $\lambda > \exp(lC/e).$
Thus, we arrive at a contradiction to our assumption that $w$ was a non-singular mixed identity. This finishes the proof.
\end{proof}

\subsection{The commutator} 
In order to deal with the many variable case under the assumptions of Theorem~\ref{thm:mixidaut}, we start with the simplest case, i.e., when $\epsilon(w)=[x,y]\coloneqq xyx^{-1}y^{-1} \in \mathbf{F}_2 = \langle x,y \rangle.$

Consider the homeomorphism $f \colon [0,1] \to [0,1]$ given by  $t.f=t^{1/2}$. Let $g_{\lambda} \colon [0,1] \to [0,1]$ be a homeomorphism with germ $t.g_{\lambda}=\lambda t$ near zero as considered above. Then
$t.[f,g_{\lambda}] = t.fg_{\lambda}f^{-1}g_{\lambda}^{-1} = \lambda x$ near zero. 

Consider now $w(x,y) = h_1 x h_2 y h_3 x^{-1} h_4 y^{-1}$. Let's assume that all constants $h_1,\ldots,h_4$ are bi-Lipschitz with derivative bounded below by $\eta>0$ near zero. We will show that $w(f,g_{\lambda})$ is not trivial for $\lambda$ large enough. We write
\begin{eqnarray*}
w(x,y) &=& [x,y] yxy^{-1}x^{-1}h_1 x h_2 y h_3 x^{-1} h_4 y^{-1} \\
&=& [x,y] yx(y^{-1}((x^{-1}h_1 x) h_2) y h_3) x^{-1} h_4 y^{-1}.
\end{eqnarray*}

We will now construct lower bounds for the derivative of $w(f,g_{\lambda})$ near zero, starting from the innermost bracket. The derivative of $f^{-1}h_1f$ is bounded below by $\eta^{1/2}$ near zero. Similarly, the derivative of $f^{-1}h_1fh_2$ is bounded below by $\eta^{1/2+1}$. Hence, the derivative of $g_{\lambda}^{-1}f^{-1}h_1fh_2g_{\lambda}$ is bounded below by $\eta^{3/2}$ near zero as well. We obtain a similar bound of $\eta^{5/2}$ for $g_{\lambda}^{-1}f^{-1}h_1fh_2g_{\lambda}h_3$ and of $\eta^5$ for $fg_{\lambda}^{-1}f^{-1}h_1fh_2g_{\lambda}h_3f^{-1}$. Finally we get the bound of $\eta^6$ for $$g_{\lambda}fg_{\lambda}^{-1}f^{-1}h_1fh_2g_{\lambda}h_3f^{-1}h_4g_{\lambda}^{-1}.$$ This implies that the derivative of $w(f,g_{\lambda})$ near zero is bounded below by $\lambda \eta^6.$ This shows that $w(f,g_{\lambda}) \neq {\rm id}_{[0,1]}$ for $\lambda$ large enough.

\subsection{Proof of Theorem~\ref{thm:mixidaut}}
Our aim is now to generalize the  technique from the previous section and to apply it to an arbitrary element $w \in {\rm Homeo}^+([0,1]) \ast \mathbf{F}_2$ with $\epsilon(w) \not \in [[\mathbf{F}_2,\mathbf{F}_2],[\mathbf{F}_2,\mathbf{F}_2]]$. Without loss of generality, we may assume that $\epsilon(w)  \in [\mathbf{F}_2,\mathbf{F}_2]$. Indeed, if $\epsilon(w) \notin [\mathbf{F}_2,\mathbf{F}_2],$ then $w(x,1)$ or $w(1,x)$ are non-singular words in one variable and the result follows from Theorem~\ref{thm:onevar}.

Let $\kappa,\lambda \in (0,\infty)$ and consider $f_\kappa, g_\lambda \in {\rm Homeo}([0,1])$ with germ $t.f_{\kappa} = t^{\kappa}$ and $t.g_\lambda=\lambda t$ near zero. Then, $[f^n_{\kappa},g^m_{\lambda}]$ for $n,m \in {\mathbb Z}$ has germ
$t \mapsto \lambda^{m(\kappa^{-n}-1)}t.$ In general, for $\mathbf{F}_2 = \langle x,y\rangle$, the derived subgroup $[\mathbf{F}_2,\mathbf{F}_2]$ is equal to the free group on the set $\{[x^n,y^m] \mid nm \neq 0\}$.

We define a group homomorphism $\alpha \colon [\mathbf{F}_2,\mathbf{F}_2] \to \mathbb R[X,X^{-1}]$
into the ring for formal real Laurent series
defined on generators by $$\alpha([x^n,y^m]) \coloneqq \log(\lambda)m(X^{-n}-1) \in \mathbb R[X,X^{-1}].$$
By the result above, the map $\alpha$ computes the logarithm of the derivative near zero of $w(f_\kappa,g_\lambda).$ For arbitrary $w$, we need to find a choice of $\kappa$ and $\lambda$ so that this logarithm of the derivative near zero does not vanish. This is clearly possible if the polynomial $\alpha(w)$ is not identically zero. Now, if $w \not \in [[\mathbf{F}_2,\mathbf{F}_2],[\mathbf{F}_2,\mathbf{F}_2]]$, then the proof of \cite{elkasapythom2014goto}*{Corollary 3.8} yields the existence of a basis of $\mathbf{F}_2$ such that $\alpha(w(f_\kappa,g_\lambda))$ is not the zero-polynomial when considered in this basis. Since a change of variable does not affect the question whether $w$ is a mixed identity, we therefore may assume that there exists $\kappa_0 >0$  so that $\alpha(w(f_{\kappa_0},g_{\lambda}))$
is a non-constant affine function in $\log(\lambda)$. 
Then, the argument in the previous section can be repeated.
This concludes the proof of Theorem~\ref{thm:mixidaut}.


In analogy to $\Aut(\mathbb{Q},<)$, we conjecture that there do not exist non-singular mixed identities for ${\rm Homeo}^+([0,1])$.


\section{Acknowledgments}
The authors would like to thank the referee for suggesting improvements and for detecting an incorrect claim about the random permutation, and for pointing us to Corollary 2.10 in~\cite{MacphersonTentSimplicity}, and Manfred Droste for pointing us to the work of Adeleke--Holland \cite{MR1269842}. 
The second and third author want to thank Henry Bradford for interesting discussions on the subject. 


\bibliographystyle{plain}
\bibliography{local}

\end{document}